\documentclass[11pt]{article}
\usepackage[a4paper,margin=20mm]{geometry}
\usepackage{mathtools}
\usepackage{authblk}

\usepackage{amsmath,amssymb}
\usepackage{tikz-cd}
\newcommand{\drpullback}{\arrow[phantom]{dr}[very near start,description]{\lrcorner}}

\makeatletter
\newcommand{\leqnomode}{\tagsleft@true\let\veqno\@@leqno}
\makeatother

\renewcommand{\section}[1]{\bigskip\noindent{\bf \thesection\ \ #1.}\stepcounter{section}}
\stepcounter{section}

\makeatletter
\renewenvironment{thebibliography}[1]
     {\bigskip\noindent{\bf \refname} \footnotesize%
      \list{\@biblabel{\@arabic\c@enumiv}}%
           {\settowidth\labelwidth{\@biblabel{#1}}%
            \leftmargin\labelwidth
            \advance\leftmargin\labelsep
            \@openbib@code
            \usecounter{enumiv}%
            \let\p@enumiv\@empty
            \renewcommand\theenumiv{\@arabic\c@enumiv}}%
      \sloppy
      \clubpenalty4000
      \@clubpenalty \clubpenalty
      \widowpenalty4000%
      \sfcode`\.\@m}
     {\def\@noitemerr
       {\@latex@warning{Empty `thebibliography' environment}}%
      \endlist}
\makeatother

\providecommand{\norm}[1]{\left| {#1}\right|}

\newcommand{\Map}{\operatorname{Map}}

\newcommand{\isopil}{\stackrel{\raisebox{0.1ex}[0ex][0ex]{\(\sim\)}}%
			{\raisebox{-0.15ex}[0.28ex]{\(\rightarrow\)}}}

\newcommand{\finiteset}{F}
\newcommand{\looop}{\gamma}
\newcommand{\base}{B}
\newcommand{\basept}{b}
\newcommand{\deloop}{\mathbf{B}}

\newcommand{\wq}{/\kern -3pt /}

\title{Free loop spaces and the Cauchy--Frobenius Lemma}

\author{Joachim Kock}
\author{Thomas Jan Mikhail}

\affil{
\vspace*{-4pt}
\small
University of Copenhagen}

\date{}

\begin{document}

\setlength{\abovedisplayskip}{9pt plus 5pt minus 3pt}
\setlength{\abovedisplayshortskip}{2pt plus 5pt minus 5pt}
\setlength\belowdisplayskip{9pt plus 4pt minus 4pt}
\setlength\belowdisplayshortskip{9pt plus 3pt minus 4pt}

\maketitle

\begin{abstract}
  We upgrade the Cauchy--Frobenius Lemma (`Burnside's Lemma')
  to a
  homotopy 
  equivalence of
  $\infty$-groupoids, essentially given by double counting/Fubini in the free
  loop space of the
  quotient.
\end{abstract}

\paragraph{Introduction.}
For a finite group $G$ acting on a finite set $\finiteset$, the 
Cauchy--Frobenius Lemma (also known as Burnside's Lemma) 
states the equation of cardinalities of finite sets
\begin{equation}\label{burnside}
\norm{\finiteset/G} = \frac{1}{\norm{G}}\sum_{g\in G} \norm{\finiteset^g} .
\end{equation}
Here $\finiteset/G$ is the {\em naive} quotient of the 
action (that is, the set of orbits), and $\finiteset^g = \{ x\in F\mid x.g=x\}$
the set of fixpoints for $g
\in G$.

That the identity is
more complicated than the formula for the homotopy 
cardinality (see~\ref{par:cardinality}) of the {\em weak} quotient 
(the action groupoid) 
\vspace*{-6pt}
\begin{equation}\label{finiteset//G}
  \norm{\finiteset\wq G} = \frac{\norm{\finiteset}}{\norm{G}} 
\end{equation}
can be explained by the interpretation $\finiteset/G = \pi_0(
\finiteset\wq G)$ and then blaming $\pi_0$ for the complication. The
appearance of $\pi_0$ is also an obstacle to realizing the identity as
the homotopy cardinality of an equivalence of groupoids.

In this note we exhibit such an equivalence (\ref{Formula 1} below), by
giving a different interpretation of $\norm{\finiteset/G}$, namely as
the homotopy cardinality of the free loop space of $F \wq G$:
$$
\norm{\finiteset/G} = \norm{\pi_0(\finiteset\wq G)} =
\norm{\Lambda (\finiteset\wq G)}.
$$

Once we focus on the free loop space, the equivalence is
straightforward to set up, amounting essentially to 
fibre sequence arguments (generalizing~\eqref{finiteset//G}),
double counting and a Fubini 
argument. In fact, the constructions work the same for 
$\infty$-groupoids (spaces), and the key insight is a general feature of
free loop spaces~\eqref{Lambda=pi0}:
\vspace*{-2pt}
$$
\norm{\Lambda X } = \norm{\pi_0 X},
$$
referring now to the homotopy cardinality of $\pi$-finite spaces, cf.~\ref{par:cardinality}.
This observation is potentially of wider applicability,
to deal with counting problems
that involve $\pi_0$.

\vspace*{-6pt}
\paragraph{Free loop spaces.} \label{sec:FreeLoopSpaces}
The {\em free loop space} $\Lambda X$ of a space $X$ is by definition the 
pullback of the diagonal along itself. The fibre over a point $x \in X$ 
is the usual based loop space:
\[
\begin{tikzcd}[row sep={3.2em,between origins}, column sep={5em,between origins}]
\Omega_x X \drpullback \ar[d] \ar[r] &
\Lambda X \drpullback \ar[d] \ar[r] &
X \ar[d, "\Delta"] \\
1 \ar[r, "x"'] &
X \ar[r, "\Delta"'] &
X\times X  .
\end{tikzcd}
\]
(The more classical definition $\Lambda X:= 
\Map(S^1,X)$ is equivalent: an element in the
pullback of $\Delta$ along itself
consists of two points $x,y\in X$ and {\em two} paths
from $x$ to $y$, forming together a loop.)

\medskip

{\em Example:} In the special case $X=\deloop G$ (the classifying space of a group 
$G$), we have
\begin{equation} \label{eq:LambdaBG}
\Lambda (\deloop G) \simeq G\wq G \simeq \sum_{g\in \operatorname{cc}(G)} \deloop
C_G(g) .
\end{equation}

\vspace*{-12pt}

\pagebreak

\noindent
The weak quotient refers to the adjoint action (action of $G$ on
itself by conjugation). The second equivalence follows from decomposition into 
connected components, where $ \operatorname{cc}(G) = \pi_0(G \wq G) $ is the set of conjugacy classes of $G$, and $ C_G(g) = \Omega_g(G \wq G) $ the centralizer.

\vspace*{-6pt}
\paragraph{Cauchy--Frobenius for $\infty$-groupoids.}
\label{par:Lambda_of_map}
The free loop space functor applied to a map $E \to B$ fits into the commutative square
\vspace*{-3pt}
\[
\begin{tikzcd}[row sep={3em,between origins}, column sep={4em,between origins}]
\Lambda E \ar[d] \ar[r] & E \ar[d]  \\
\Lambda \base \ar[r] & \base
\end{tikzcd}
\]
(not generally a pullback).
We can now write $\Lambda E$ as a homotopy sum of its fibres in three 
ways: over $E$, over $\Lambda \base$, or over $\base$. 
Homotopy double counting with respect to $E \leftarrow \Lambda 
E \to \Lambda \base$ gives the equivalences
\vspace*{-4pt}
\begin{equation}\leqnomode\label{Formula 0}\tag{\sc Formula 0}
\int_{x\in E} 
\Omega_x E
\simeq
\Lambda E 
\simeq
\int_{\looop  \in \Lambda \base} F^\looop   .
\end{equation}
Here by definition, $F^\looop$ is the `fixpoint space' for a loop 
$\looop\in \Lambda \base$ (based at a point $\basept \in \base$), 
formally defined together with the `incidence 
correspondence' $\mathbb{I}_\basept$ as the pullbacks
\begin{equation}\label{defn:I}
\begin{tikzcd}[row sep={3em,between origins}, column sep={4em,between origins}]
F^\looop   \drpullback \ar[d] \ar[r] &
\mathbb{I}_\basept \drpullback \ar[d] \ar[r] & 
\Lambda E \ar[d] \\
1 \ar[r, "\looop  "'] &
\Omega_\basept \base \drpullback \ar[d] \ar[r]
& \Lambda \base \ar[d] \\
& 1 \ar[r, "\basept"'] & \base ,
\end{tikzcd}
\end{equation}

\vspace*{-6pt}

\noindent further analyzed in~\ref{fixpt} below.

Continuing with the right-hand side of~\ref{Formula 0}, since the domain 
of integration $\Lambda\base$ splits as the homotopy sum of
its fibres over $\base$, we can now apply Fubini:
$$
\int_{\looop  \in \Lambda \base} F^\looop   \simeq \int_{\basept\in \base}  
\int_{\looop  \in \Omega_\basept \base} F^\looop   .
$$
Altogether we get the {\em Cauchy--Frobenius Lemma for 
$\infty$-groupoids}, expressing two perspectives on $\Lambda E$: 
\begin{equation}\leqnomode\label{Formula 1}\tag{\sc Formula 1}
\boxed{
\int_{x\in E} 
\Omega_x E
\simeq
\Lambda E 
\simeq
\int_{\basept \in \base}  \int_{\looop  \in \Omega_\basept \base} F^\looop  
}
\end{equation}

\paragraph{Fixpoint interpretation.} \label{fixpt}
Let $(\base, \ast)$ be a pointed connected space. In~\eqref{defn:I} we defined $\mathbb{I}$
as the pullback below left:
\vspace*{-2pt}
\[
\begin{tikzcd}[row sep={3em,between origins}, column sep={4em,between origins}]
\mathbb{I} \drpullback \ar[d] \ar[r] & 
\Lambda E \ar[d] \\
\Omega \base \drpullback \ar[d] \ar[r]
& \Lambda \base \ar[d] \\
1 \ar[r, "\ast"'] & \base 
\end{tikzcd}
\qquad = \qquad
\begin{tikzcd}[row sep={3em,between origins}, column sep={4em,between origins}]
\mathbb{I} \drpullback \ar[d] \ar[r] & 
\Lambda E \ar[d] \\
F \drpullback \ar[d] \ar[r]
& E \ar[d] \\
1 \ar[r, "\ast"'] & \base 
\end{tikzcd}
\]
and consequently we could also have defined it by the pullbacks above 
right, where $F$ denotes the fibre of $E\to\base$.
Composing with the pullback that defines the free loop space, we thus have the 
pullback characterization of $\mathbb{I}$ as below left:
\begin{equation}\label{shear}
\begin{tikzcd}[row sep={3em,between origins}, column sep={5em,between origins}]
\mathbb{I} \ar[d] \ar[r] \drpullback & \Lambda E \drpullback  
\ar[d] \ar[r] & E \ar[d, "\Delta"]  \\
F \ar[r] & E \ar[r, "\Delta"'] & E\times E
\end{tikzcd}
\qquad = \qquad
\begin{tikzcd}[row sep={3em,between origins}, column sep={5em,between origins}]
\mathbb{I} \ar[d] \ar[r] \drpullback & F\times_E F \drpullback  
\ar[d] \ar[r] & E \ar[d, "\Delta"]  \\
F \ar[r, "\Delta"'] & F\times F \ar[r] & E\times E
\end{tikzcd}
\end{equation}
and consequently a pullback characterization of $\mathbb{I}$ as above 
right.

\pagebreak

On the other hand, the cube of pullback squares
\[
\begin{tikzcd}[column sep={2.9em,between origins}, row sep={2em,between origins}]
	& F && {F\times_E F} \\
	E && F \\
	& 1 && {\Omega\base} \\
	\base && 1
	\arrow[from=1-2, to=2-1]
	\arrow[from=1-2, to=3-2]
	\arrow[from=1-4, to=1-2]
	\arrow[from=1-4, to=2-3]
	\arrow[from=1-4, to=3-4]
	\arrow[from=2-1, to=4-1]
	\arrow[from=2-3, to=2-1]
	\arrow[from=2-3, to=4-3]
	\arrow[from=3-2, to=4-1]
	\arrow[from=3-4, to=3-2]
	\arrow[from=3-4, to=4-3]
	\arrow[from=4-3, to=4-1]
\end{tikzcd}
\]
gives a natural identification $F\times_E F \simeq F \times 
\Omega \base$, 
under which the projections become the shear map
\begin{equation*}
  F\times \Omega \base \longrightarrow F \times F, \qquad
  (x,\looop) \longmapsto (x, x.\looop) .
\end{equation*}
Coming back to~\eqref{shear}, we see that $\mathbb{I}$ is
the pullback of the shear map along the diagonal, and hence is the 
incidence correspondence, 
which is the space of $(x,\looop,p)$ with $p:x.\looop\isopil x$.
Its fibre $F^\looop$ is the space of pairs $(x,p)$ with $p : x.\looop 
\isopil x$, the space of (homotopy) fixpoints of $\gamma \in \Omega B$.

\vspace*{-6pt}
\paragraph{Homotopy cardinality.} \label{par:cardinality}
A space is called {\em $\pi$-finite} if it is truncated, has finitely many connected components, and all its homotopy groups are finite. The {\em homotopy cardinality} of $\pi$-finite spaces 
is characterized 
by the recursive identity
\vspace*{-4pt}
$$
\norm{X} = \sum_{x\in \pi_0 X} \frac{1}{\norm{\Omega_x X}} 
$$
and the initial condition $\norm{1}=1$. (In the discrete case, homotopy cardinality 
agrees with usual cardinality of a finite set.)
More generally, given a map of $\pi$-finite spaces $Y\to 
X$,
we have
\begin{equation}\label{Y}
\norm{Y} = \norm{\int_{x\in X} Y_x} = \sum_{x\in \pi_0 X} \frac{\norm{Y_x}}{\norm{\Omega_x X}} .
\end{equation}

Coming back to the left-hand side of~\ref{Formula 1}, referring to the 
map $\Lambda E \to E$ (with $E$ $\pi$-finite), we get the fundamental identity
\begin{equation}\label{Lambda=pi0}
\norm{\Lambda E} = \norm{\int_{x \in E} \Omega_x E} = \sum_{x\in \pi_0 E} \frac{\norm{\Omega_x E}}{\norm{\Omega_x E}}
= \sum_{x\in \pi_0 E} 1 = \norm{\pi_0 E} .
\end{equation}

\vspace*{-8pt}
\paragraph{The classical case.} \label{classical}
We take cardinality of~\ref{Formula 1} in the special case where $\base=\deloop G$, the classifying space of a finite group $G$.
Let $G$ act on a finite
set $F$ and let $E:=F\wq G$ be the weak quotient (so that $\pi_0 E=
F/G$); the associated map $F \wq G \to \deloop G$ has fibre $F$. 
For the left-hand side of~\ref{Formula 1} we already calculated in~\eqref{Lambda=pi0}
that $\norm{\Lambda (F\wq G)} = \norm{\pi_0 (F\wq G)} = \norm{F/G}$. On the right-hand
side of~\ref{Formula 1} we apply~\eqref{Y} twice: 
the 
integral over $B = \deloop G$
becomes division by $\norm{G}$, and the second integral becomes a discrete sum over
$\Omega B \simeq G$. 
We thus arrive at
$$
  \norm{F / G} = \frac{1}{\norm{G}} \sum_{g \in G} \norm{F^g}.
$$

It remains to justify that this $F^g$
is indeed the usual fixpoint set as in
\eqref{burnside}. 
Since
$F$ and $G$ are discrete, also $F\times_E F \simeq F\times G$ is
discrete. It follows from the pullback description in~\eqref{shear}
that also $\mathbb{I}$ is discrete: we have $\mathbb{I} \simeq \{(x,g)
\in F\times G \, | \, x.g = x \} $.
So indeed we have $F^g \simeq \{x \in F \, |
\, x.g = x \} $, again discrete.
Altogether, we get the Cauchy--Frobenius Lemma~\eqref{burnside}.

(Note also that the homotopy
cardinality of~\ref{Formula 0} gives the following variant of 
Cauchy--Frobenius:
$
\norm{\finiteset/G} =
\sum_{g\in \operatorname{cc}(G)} \norm{\finiteset^g}/\norm{C_G(g)}.
$
For the trivial action of $G$ on a point, this gives the class equation
$1 = \sum_{g \in \operatorname{cc}(G)} 
\norm{C_G(g)}^{-1}$ (which is also the homotopy cardinality of 
\eqref{eq:LambdaBG}).)

\vspace*{-6pt}
\paragraph{Notes.}
Burnside~\cite{Burnside:1897} stated the lemma with attribution to 
Frobenius~\cite{Frobenius:1887}, but that 
reference was left out in the second edition of his book, and the 
result became known as Burnside's Lemma. According
to Neumann~\cite{Neumann:1979},
a version of the lemma had been published decades earlier by 
Cauchy~\cite{Cauchy:1845}.

\pagebreak

Homotopy cardinality was advocated by Joyal and by
Baez--Dolan~\cite{Baez-Dolan:finset-feynman} as a starting point for 
homotopy combinatorics,
but the 
definition goes back at least to 
Quinn~\cite{Quinn}. We refer to
\cite{Galvez-Kock-Tonks:1602.05082} for the few results required here.

Double counting 
was held by 
Aigner~\cite{Aigner:reprint} 
as one of the most important principles of enumerative combinatorics.
Homotopy double counting and the companion Fubini lemma were advocated by 
G\'alvez--Kock--Tonks~\cite[Lemmas 3.8 and 3.9]{Galvez-Kock-Tonks:1207.6404}.
The resulting identities of rational numbers sometimes exhibit
striking 
cancellations of symmetry factors, as illustrated in 
\cite{Galvez-Kock-Tonks:1207.6404} with the Fa\`a di Bruno formula for 
free operads.

The bundle viewpoint on (higher) group actions is standard in 
homotopy theory; we refer in particular to 
\cite{Nikolaus-Schreiber-Stevenson:1207.0248} for the equivalence
$F\times_E F \simeq F \times \Omega \base$
and the 
homotopy-fixpoint interpretation.

The brevity of this note exempts us from any survey of the great 
importance of free loop spaces in modern homotopy theory, for example 
in string topology and K-theory. We must 
mention, however, that
in geometry, free loop spaces appear in particular in the guise of
inertia orbifolds (where the notation $\Lambda$ originates).
Analogues of
\ref{Formula 0} and~\ref{Formula 1} are known in that
context~\cite{Lupercio-Uribe:0110207}.

Constructions similar to free
loop spaces exist in combinatorics, with permutations playing the role
of loops. One instance is Joyal's tilde construction on a combinatorial
species~\cite{Joyal:1981}. He used `Burnside's Lemma' to establish its
key property $|\widetilde M|=\norm{\pi_0 M}$ (which corresponds to the
insight~\eqref{Lambda=pi0} of our proof).
Exploiting related viewpoints in probability theory,
Baez~\cite{Baez:2412.16386} gives a (groupoid-level) formula for the
expectation of a set-valued functor on $G\wq G$. We point out that 
his 
formula specializes to the Cauchy--Frobenius Lemma when the functor is
the fixpoint statistics of a group action, and that 
the
formula follows
from (the right-hand part of) \ref{Formula 1} if only $\Lambda
E\to\Lambda \base$ is replaced by a general map $X\to\Lambda \base$
(same proof).

\vfil
\vfill

\footnotesize

\noindent 
\hrule
\smallskip
\noindent Funding:
Both authors were funded by grant No.~10.46540/3103-00099B from the Independent Research
  Fund Denmark. We also acknowledge support from research grant PID2020-116481GB-I00
  (AEI/FEDER, UE) of Spain and grant 2021-SGR-1015 of Catalonia, and through the
  Severo Ochoa and Mar\'ia de Maeztu Program for Centers and Units of Excellence in
  R\&D grant number CEX2020-001084-M as well as the Danish National Research
  Foundation through the Copenhagen Centre for Geometry and Topology (DNRF151).

\end{document}